\documentclass[english]{article}
\usepackage[T1]{fontenc}
\usepackage[latin1]{inputenc}
\usepackage{verbatim}
\usepackage{amsmath}
\usepackage{amssymb}

\makeatletter
\usepackage{amsfonts}

\providecommand{\U}[1]{\protect\rule{.1in}{.1in}}
\newtheorem{theorem}{Theorem}

\newtheorem{definition}[theorem]{Definition}

\newtheorem{lemma}[theorem]{Lemma}

\newtheorem{proposition}[theorem]{Proposition}
\newtheorem{remark}[theorem]{Remark}

\newenvironment{proof}[1][Proof]{\noindent\textbf{#1.} }{\ \rule{0.5em}{0.5em}}

\usepackage{babel}
\makeatother
\begin{document}

\title{The Poisson Kernel for Hardy Algebras}

\author{Paul S. Muhly%
\thanks{Supported in part by grants from the National Science Foundation and
from the U.S.-Israel Binational Science Foundation.%
}\\
Department of Mathematics\\
University of Iowa\\
Iowa City, IA 52242\\
e-mail: muhly@math.uiowa.edu \and Baruch Solel%
\thanks{Supported in part by the U.S.-Israel Binational Science Foundation
 and by the B. and G. Greenberg Research Fund (Ottawa).%
}\\
Department of Mathematics\\
Technion\\
32000 Haifa, Israel\\
e-mail: mabaruch@techunix.technion.ac.il}

\maketitle

\section{Introduction}

This note contributes to a circle of ideas that we have been developing
recently in which we view certain abstract operator algebras, which
we call Hardy algebras, and which are noncommutative generalizations
of classical $H^{\infty}$, as spaces of functions defined on their
spaces of representations \cite{MS04, MS05, MS05a, MS06}. This perspective leads to a number of pleasant
formulas that are very reminiscent of formulas from complex function
theory on the unit disc. More important, however, they help to reveal
structural properties of the algebras and they help to clarify the
interplay among various constructs that are at work in their analysis.
Even in the classical setting of complex functions of one variable,
insight is sometimes gained by viewing classical $H^{\infty}$ as
a space of functions on its space of representations, which are parameterized,
essentially, by all the completely non-unitary contractions. Another
source of motivation is the work of Popescu, Davidson and Pitts, and
others who have done extensive work on free semigoup algebras.\footnote{For a nice survey of the basics of free semigroup algebras, we recommend Ken Davidson's article \cite{kD01}.} Indeed,
many of the results that we prove here have been been anticipated
in this work. What is novel about our approach, however, is the systematic
use of ``duality of correspondences'' to put into evidence the effectiveness
of viewing elements of our Hardy algebras as functions on operator
discs. When this is done, proofs in the free semigroup picture often
become simpler, shorter and more perspicuous. And they extend to a wide variety of additional
situations in the literature that are of interest.

In the next section, we introduce the basic players in our theory: a $W^*$-algebra $M$, a $W^*$-correspondence $E$ over $M$ and the Hardy algebra they generate: $H^{\infty}(E)$.  We then describe how a normal representation $\sigma:M\to B(H_{\sigma})$ gives rise to a ``dual'' correspondence, denoted $E^{\sigma}$ and we describe how elements of $H^{\infty}(E)$ may be realized as functions defined on the unit ball of the space of adjoints of $E^{\sigma}$, $\mathbb{D}((E^{\sigma})^*)$. In Section 3, we define a generalization of the Poisson kernel, which ``reproduces'' the values on $\mathbb{D}((E^{\sigma})^*)$ of the ``functions'' coming from $H^{\infty}(E)$.  When $M = E = \mathbb{C}$ and $\sigma$ is the one dimensional representation of $M$, then $H^{\infty}(E)$ is classical $H^{\infty}$ realized as analytic Toeplitz operators, and our Poisson kernel is easily seen to be the classical Poison kernel formulated in terms of operators on Hilbert space.  Our representation theorems, Theorems \ref{Reproducing1} and \ref{Reproducing2} are easily seen to be natural generalizations of the Poisson integral formuala.  They also are easily seen to be generalizations of formulas that Popescu developed in \cite{gP99} and elsewhere, and they are closely related to formulas that Arveson developed in \cite{wA00}.  In the fourth section, we relate our Poisson kernel to the idea of a characteristic operator function and show how the Poisson kernel identifies the ``model space'' for the canonical model that can be attached to a point in the disc $\mathbb{D}((E^{\sigma})^*)$ - a structure we developed in \cite{MS05}.  We were inspired here by \cite{gP06} and other results from literature. In the next section, Section 5, we consider a Poisson kernel on the unit ball of $E$, $\mathbb{D}(E)$.  Owing to our duality theorem \cite[Theorem 3.9]{MS04}, one can think of this ball as the place to evaluate elements in $H^{\infty}(E^{\sigma})$, but in addition, it captures ideas about ``(left) point evaluations'' that appear in the systems theory literature, cf. \cite{ADD90}.  Finally, in Section 6, we connect our Poission kernel to the idea of curvature and complement results that we proved in \cite{MS03}.  Again, our analysis extends parts of the theory of curvature for not-necessarily-commuting row contractions that was developed by Popescu in \cite{gP01}. His work, in turn, was based on investigations by Arveson \cite{wA00} in which he introduced a notion of curvature to study properties of commuting row contractions.

\section{Preliminaries}

We recall a few key ideas from \cite{MS04} and we refer to that paper
for further discussion and references about the setup with which we
will be working here. Throughout this note $M$ will be a fixed $W^{*}$-algebra.
We also fix a $W^{*}$-correspon\-dence $E$ over $M$. This means that
$E$ is a self-dual Hilbert $C^{*}$-module over $M$ and that there
is a normal homomorphism $\varphi$ from $M$ into the $W^{*}$-algebra
of all continuous module maps on $E$, $\mathcal{L}(E)$, giving $E$
an action of $M$ that makes $E$ a bimodule over $M$. We shall form
the (balanced) tensor powers of $E$, $E^{\otimes n}$, which are
all $W^{*}$-correspondences over $M$, and we shall denote the left
action of $M$ on $E^{\otimes n}$ by $\varphi_{n}$. It is defined
by the formula\[
\varphi_{n}(a)(\xi_{1}\otimes\xi_{2}\otimes\cdots\otimes\xi_{n})=(\varphi(a)\xi_{1})\otimes{(\xi}_{2}\otimes\cdots\otimes\xi_{n}).\]
We shall write $E^{\otimes0}=M$, viewed as a bimodule over itself,
so in particular, $\varphi_{0}(a)\xi=a\xi$. The direct sum $E^{\otimes0}\oplus E^{\otimes1}\oplus E^{\otimes2}\oplus\cdots$
is a $W^{*}$-correspondence over $M$ in an obvious and natural way,
which we shall denote by $\mathcal{F}(E)$ and call the \emph{Fock
space} over $E$. The left action of $M$ on $\mathcal{F}(E)$ is
the sum of the $\varphi_{n}$ and will be denoted $\varphi_{\infty}$.
Thus, for $a\in M$,\[
\varphi_{\infty}(a)=\textrm{diag}(\varphi_{0}(a),\varphi_{1}(a),\varphi_{2}(a),\ldots),\]
when we view operators as matrices on $\mathcal{F}(E)$ as we shall.
An element $\xi\in E$ defines a \emph{creation operator} $T_{\xi}$
on $\mathcal{F}(E)$ via the formula $T_{\xi}\eta=\xi\otimes\eta$.
This operator is bounded, with adjoint given by the formula $T_{\xi}^{*}(\zeta\otimes\eta)=\varphi_{\infty}(\langle\xi,\zeta\rangle)\eta$.
Matricially, $T_{\xi}$ has a form of an operator-valued weighted
shift:\[
T_{\xi}=\left[\begin{array}{ccccc}
0\\
T_{\xi}^{(1)} & 0 &  & 0\\
0 & T_{\xi}^{(2)} & 0\\
 & 0 & T_{\xi}^{(3)} & \ddots\\
0 &  & \ddots & \ddots & \ddots\end{array}\right]\]
 where $T_{\xi}^{(n)}$ maps $E^{\otimes(n-1)}$ into $E^{\otimes n}$
by tensoring with $\xi$. The ultraweakly closed subalgebra of $\mathcal{L}(\mathcal{F}(E))$
generated by the $T_{\xi}$, $\xi\in E$, and the $\varphi_{\infty}(a)$,
$a\in M$, is called the \emph{Hardy algebra} of $E$ and is denoted
$H^{\infty}(E)$. Numerous examples of Hardy algebras may be found
in the literature that we cite, and elsewhere, so we won't go into
detail here. However, we do want to point out that when $M=E=\mathbb{C}$,
the complex numbers, then $H^{\infty}(E)$ is the classical Hardy
space of bounded analytic functions on the open unit disc, $H^{\infty}$,
realized as the algebra of all (bounded) analytic Toeplitz operators
on the space $\ell^{2}(\mathbb{Z}_{+})$. Hence the terminology.

A fundamental feature of our theory is that the ultraweakly continuous
completely contractive representations of $H^{\infty}(E)$ can be
parametrized by the normal representations of $M$ and certain
contraction operators in a fashion that we want to describe in some
detail. Let $\sigma:M\to B(H)$ be a normal representation
of $M$ on a Hilbert space $H$. Then $\sigma$ induces a normal representation
$\sigma^{E}$ of $\mathcal{L}(E)$ on $E\otimes_{\sigma}H$, defined
via the formula $\sigma^{E}(X)=X\otimes I_{H}$. In fact, $\sigma^{E}$
\emph{is} called the induced representation of $\mathcal{L}(E)$ determined
by $\sigma$, and we refer to \cite{mR74} for a discussion of the
general theory. If we form $\sigma^{E}\circ\varphi$ we obtain a new
representation of $M$ that we denote simply by $\varphi\otimes I$
and refer to as \emph{the induced representation of $M$ determined
by} $\sigma$ (and $E$). Suppose that $T$ is an operator from $E\otimes_{\sigma}H$
to $H$ of norm at most one that intertwines the induced representation
of $M$ and $\sigma$, i.e., suppose \begin{equation}
T(\varphi(a)\otimes I_{H})=\sigma(a)T\label{eq:intertwine}\end{equation}
for all $a\in M$, then $T$ determines an ultraweakly continuous,
completely contractive bimodule map $\hat{T}$ from $E$ to $B(H)$
via the formula\begin{equation}
\hat{T}(\xi)h=T(\xi\otimes h),\label{eq:hat}\end{equation}
$\xi\in E$ and $h\in H$. That is, $\hat{T}:E\to B(H)$ is completely
contractive, where $E$ is regarded as an operator space in the operator
space structure it inherits as a subspace of its linking algebra \cite[Page 398]{MS98},
and is continuous with respect to the natural so-called $\sigma$-topology
of \cite{BBLS04} and the ultraweak topology on $B(H)$. The bimodule
property refers to the equation $\hat{T}(\varphi(a)\xi b)=\sigma(a)\hat{T}(\xi)\sigma(b)$,
which is satisfied for all $a,b\in M$ and $\xi\in E$. We call the
pair $(\hat{T},\sigma)$ an \emph{(ultraweakly continuous completely
contractive) covariant representation} of $E$ (and $M$) on $H$.
Conversely, given such a representation of $E$ and $M$ on a Hilbert
space $H$, $(S,\sigma)$, the formula \begin{equation}
\tilde{S}(\xi\otimes h):=S(\xi)h,\label{eq:tilde}\end{equation}
$\xi\otimes h\in E\otimes_{\sigma}H$ defines an operator of norm
at most $1$ from $E\otimes_{\sigma}H$ to $H$ that satisfies equation
(\ref{eq:intertwine}). We denote this operator by $\tilde{S}$, i.e.,
$\tilde{S}(\varphi(a)\otimes I_{H})=\sigma(a)\tilde{S}$ for all $a\in M$.
Clearly, we have $\tilde{\hat{T}}=T$ and $\hat{\tilde{T}}=T$.

The key point is that each ultraweakly continuous, completely contractive
representation $\rho$, say, of $H^{\infty}(E)$ on a Hilbert space
$H$ determines a completely contractive covariant representation
of $E$ and $M$ on $H$ through the formulas\[
\sigma(a)=\rho(\varphi_{\infty}(a))\]
and\[
T(\xi)=\rho(T_{\xi}),\]
and conversely, (almost) every completely contractive covariant representation
$(T,\sigma)$ ``integrates'' to an ultraweakly continuous, completely
contractive representation $\rho$ through these formulas. We say
``almost'' because while every $(T,\sigma)$ ``integrates'' to a \emph{norm}-continuous,
completely contractive representation $\rho$ of the \emph{norm}-closed
algebra generated by $\{ T_{\xi}\}_{\xi\in E}$ and $\varphi_{\infty}(M)$,
which we denote by $\mathcal{T}_{+}(E)$ and call the \emph{tensor
algebra} of the correspondence, the representation $\rho$ need not
extend all the way to $H^{\infty}(E)$. (We will say more about this
in a moment.) We write $\sigma\times T$ for the representation determined
by $(T,\sigma)$ on the norm-closed algebra whether or not it extends
to $H^{\infty}(E)$%
\footnote{In some of our papers, we have written ``$T\times\sigma''$instead
of ``$\sigma\times T$''. We apologize for any confusion this may
create.%
}. If $\Vert\tilde{T}\Vert<1,$ then $\sigma\times T$ \emph{does}
extend to an ultraweakly continuous, completely contractive representation
of $H^{\infty}$ \cite[Corollary 2.14]{MS04}. Thus we can say that
once a normal representation $\sigma$ of $M$ on $H$ is
given, then there is a bijective correspondence between the strictly
contractive intertwiners of $\sigma$ and $\varphi\otimes I$ and
the ultraweakly continuous, completely contractive representations
$\rho$ of $H^{\infty}(E)$ on $H$ such that $\rho\circ(\varphi_{\infty}\otimes I_{H})=\sigma$
and such that $\Vert\rho(T_{\xi})\Vert<c\Vert\xi\Vert$ for all $\xi\in E$,
where $c$ is a prescribed constant less than $1$. This observation
suggests that we may adopt the perspective of viewing elements of
$H^{\infty}(E)$ as functions on the space of (ultraweakly continuous,
completely contractive) representations of $H^{\infty}(E)$ in a concrete
and transparent fashion. This suggestion was the principal point of
\cite{MS04} and has been the focus of much of our subsequent work.
To help explain further the functional perspective initiated in \cite{MS04},
we require the following definition.

\begin{definition}\label{Esigma}If $\sigma:M\to B(H)$ is a normal
 representation of $M$ on the Hilbert space $H$, then we
define $E^{\sigma}$ to be the space of bounded operators $\eta:H\to E\otimes_{\sigma}H$
with the property that $\eta\sigma(a)=(\varphi(a)\otimes I_{H})\eta$
for all $a\in M$. We call $E^{\sigma}$ the $\sigma$\emph{-dual}
of $E$. We write $\mathbb{D}(E^{\sigma})$ for the \emph{open} unit
ball in $E^{\sigma}$.\end{definition}

Evidently, the elements of $E^{\sigma}$are precisely the adjoints
of the space of operators that satisfy equation (\ref{eq:intertwine}).
Suppose $\eta\in\mathbb{D}(E^{\sigma})$ is given. Then $\eta^{*}$
satisfies equation (\ref{eq:intertwine}) and determines an ultraweakly,
completely contractive covariant representation $(\widehat{\eta^{*}},\sigma)$
of $E$ on $H$. Further, with the aid of \cite[Corollary 2.14]{MS04},
the formulas\[
\varphi_{\infty}(a)\to\sigma(a)\]
and\[
T_{\xi}\to\widehat{\eta^{*}}(\xi)\]
extend to give an ultraweakly continuous, completely contractive representation
$\sigma\times\widehat{\eta^{*}}$ of $H^{\infty}(E)$ on $H$. On
elements of the form ${(T}_{\xi_{1}}\otimes I_{H}){(T}_{\xi_{2}}\otimes I_{H})\cdots{(T}_{\xi_{n}}\otimes I_{H})={(T}_{\xi_{1}\otimes\xi_{2}\otimes\cdots\xi_{n}}\otimes I_{H})$,
for example, $\sigma\times\widehat{\eta^{*}}$ is given by the formula\begin{equation}
\sigma\times\widehat{\eta^{*}}{(T}_{\xi_{1}\otimes\xi_{2}\otimes\cdots\xi_{n}}\otimes I_{H})=\widehat{\eta^{*}}(\xi_{1})\widehat{\eta^{*}}(\xi_{2})\cdots\widehat{\eta^{*}}(\xi_{n}).\label{eq:evaluation}\end{equation}
Following \cite{MS05}, we introduce the following terminology.

\begin{definition}For $\eta\in\mathbb{D}(E^{\sigma})$ and for $X\in H^{\infty}(E)$,
we define\begin{equation}
\widehat{X}(\eta^{*}):=\sigma\times\widehat{\eta^{*}}(X).\label{eq:Xhat}\end{equation}
The resulting function $\widehat{X}:\mathbb{D}(E^{\sigma})^{*}\to B(H)$
is called the \emph{Fourier transform} of $X$.\end{definition}

Perhaps the term ``$Z$-transform'' is preferable to ``Fourier transform'',
but both conjure up formulas such as $\widehat{XY}(\eta^{*})=\widehat{X}(\eta^{*})\widehat{Y}(\eta^{*})$
that are clearly evident from (\ref{eq:evaluation}).

\begin{remark}Suppose $M=E=\mathbb{C}$ and that $H$ also is $\mathbb{C}$.
Then of course $\sigma$ can only be the identity representation of
$M=\mathbb{C}$ on $H$, $E^{\sigma}$ also may be identified with
$\mathbb{C}$. In this situation, then, $\mathbb{D}(E^{\sigma})$
is just the open unit disc $\mathbb{D}$ in the complex plane. The
Fourier transform takes an $X$ in $H^{\infty}(E)$, which by our
definition is an infinite, lower-triangular, Toeplitz matrix on $\ell^{2}(\mathbb{Z}_{+})$\[
X=\left(\begin{array}{ccccc}
a_{0} & 0 & 0 & \dots & \dots\\
a_{1} & a_{0} & 0 & \ddots & \ddots\\
a_{2} & a_{1} & a_{0} & 0 & \ddots\\
a_{3} & a_{2} & a_{1} & \ddots & \ddots\\
\vdots & \ddots & \ddots & \ddots & \ddots\end{array}\right)\]
that represents a bounded operator, to a function from $\mathbb{D}$
to operators on $H=\mathbb{C}$, i.e., to numbers. To compute them,
simply note that for $\eta\in\mathbb{D}$, $\eta^{*}$ is just the
complex conjugate of $\eta$, $\overline{\eta}$, and equation (\ref{eq:Xhat})
implies that $\widehat{X}(\eta^{*})$ is nothing but multiplication
by the complex number $\sum_{k=0}^{\infty}a_{k}\overline{\eta}^{k}$
on $\mathbb{C}$, i.e., for $c\in\mathbb{C}$, $\widehat{X}(\eta^{*})c=(\sum_{k=0}^{\infty}a_{k}\overline{\eta}^{k})c.$
It is clear in this example, that for no $\eta$ on the boundary of
$\mathbb{D}$ does $\sigma\times\widehat{\eta^{*}}$ extend to an
\emph{ultraweakly continuous} representation of $H^{\infty}(E)$.
If, next, $H=\mathbb{C}^{n}$, and again if $\sigma(a)\xi=a\xi$,
for $a\in M=\mathbb{C}$, then $E^{\sigma}$ may be viewed as the
$n\times n$ matrices over $\mathbb{C}$, and $\mathbb{D}(E^{\sigma})$
consists of all those $n\times n$ matrices of norm less than $1$.
If $T$ is such a matrix, then $\widehat{X}(T^{*})$ is the operator
on $H=\mathbb{C}^{n}$ given by a similar formula: \begin{equation}
\widehat{X}(T^{*})\xi=(\sum_{k=0}^{\infty}a_{k}T^{*^{k}})\xi.\label{eq:FunctionEval}\end{equation}
 It is clear in this case, that for $\Vert T\Vert=1$, $\sigma\times\widehat{T^{*}}$
extends to an ultraweakly continuous representation of $H^{\infty}(E)$
on $H$ if and only if the spectral radius of $T$ is less than one.
Finally, if $H$ is an infinite dimensional Hilbert space, so that
$\sigma(a)\xi=a\xi$, as before, then $E^{\sigma}$ may be identified
with $B(H)$ and $\mathbb{D}(E^{\sigma})$ may be viewed as the collection
of all operators on $H$ of norm less than one. In this case, $\widehat{X}(T^{*})$
again is given by the formula (\ref{eq:FunctionEval}). Now, however,
the $T$'s of norm one for which $\sigma\times\widehat{T^{*}}$ extends
to an ultraweakly continuous representation of $H^{\infty}(E)$ are
precisely those whose minimal unitary dilations are absolutely continuous
with respect to Lebesgue measure on the circle. Such a contraction
splits into a completely non-unitary contraction and an absolutely
continuous unitary operator. The value of $\widehat{X}(T^{*})$ for
such a $T$ is given by the Sz.-Nagy - Foia\c{s} functional calculus.
In \cite[Section 7]{MS04} we showed, in general, that if $\eta\in\mathbb{D}(E^{\sigma})$
is such that $\eta^{*}$ is ``completely noncoisometric'', then $\sigma\times\widehat{\eta^{*}}$
extends to an ultraweakly continuous representation of $H^{\infty}(E)$.
Beyond this, it is a mystery to us about how to identify points $\eta$
on the boundary of $\mathbb{D}(E^{\sigma})$ in general such that
$\sigma\times\widehat{\eta^{*}}$ extends to an ultraweakly continuous
representation of $H^{\infty}(E)$.\end{remark}

The reason we focus on $E^{\sigma}$ rather than on the space of its
adjoints, $(E^{\sigma})^{*}$, at least for some purposes, is that
$E^{\sigma}$ is a $W^{*}$-correspondence over the commutant of $\sigma(M),$
$\sigma(M)'$. The point to keep in mind is that the commutant of
$\sigma^{E}(\mathcal{L}(E))$ is $I_{E}\otimes\sigma(M)'$ \cite[Theorem 6.23]{mR74},
and so $E^{\sigma}$ becomes a bimodule over $\sigma(M)'$ according
to the formula\[
a\cdot\eta\cdot b:=(I_{E}\otimes a)\eta b,\]
$a,b\in\sigma(M)'$, $\xi\in E^{\sigma}$. The $\sigma(M)'$-valued
inner product on $E^{\sigma}$ is given simply by operator multiplication:\[
\langle\eta,\zeta\rangle:=\eta^{*}\zeta,\]
$\eta,\zeta\in E^{\sigma}$. For more details about the structure
of $E^{\sigma}$and examples, see Sections 3 and 4 of \cite{MS04}.

One of the important points for us in this note is that for the representations
$\rho$ of $H^{\infty}(E)$ that we defined in \cite{MS99} and called
\emph{induced representations}, the commutant of $\rho(H^{\infty}(E))$
can be expressed in terms of induced representations of $H^{\infty}(E^{\sigma})$\emph{.}

\begin{definition}\label{InducedCovariant}Let $\sigma:M\to B(H)$
be a normal representation of $M$ on a Hilbert space $H$
and form the Hilbert space $\mathcal{F}(E)\otimes_{\sigma}H$. The
\emph{induced covariant representation} of $E$ determined by $\sigma$
is the representation $(V,\varphi_{\infty}\otimes I_{H})$ where $V:E\to B(\mathcal{F}(E)\otimes_{\sigma}H)$
is defined by the equation\[
V(\xi)(\eta\otimes h):=(\xi\otimes\eta)\otimes h,\]
$\xi\in E$, and $\eta\otimes h\in\mathcal{F}(E)\otimes_{\sigma}H$.
The integrated form of $(V,\varphi_{\infty}\otimes I_{H})$, $(\varphi_{\infty}\otimes I_{H})\times V$,
is called the representation of $H^{\infty}(E)$ \emph{induced} by
$\sigma$. We shall usually write $\sigma^{\mathcal{F}(E)}$ for $(\varphi_{\infty}\otimes I_{H})\times V$,
and most frequently, we will simply write $X\otimes I_{H}$ for $\sigma^{\mathcal{F}(E)}(X)$,
$X\in H^{\infty}(E)$.\end{definition}

The map $V$ is essentially the map defining the tensor powers of
$E$ and the associated map $\tilde{V}:E\otimes(\mathcal{F}(E)\otimes H)\to\mathcal{F}(E)\otimes H$
appears to be just the identity map embedding $\sum_{k=1}^{\infty}{(E}^{\otimes k}\otimes_{\sigma}H)$
into $\mathcal{F}(E)\otimes_{\sigma}H$. However, it is a bit more
complicated. There is a shift involved, as we shall see later in equation
(\ref{eq:Shift}) and subsequent analysis.

Observe that if $\eta\in E^{\sigma}$, then for each $k\geq0$, $I_{E^{\otimes k}}\otimes\eta$
may be viewed as a map from $E^{\otimes k}\otimes_{\sigma}H$ to $E^{\otimes(k+1)}\otimes_{\sigma}H$.
Further, due to the balanced nature of the tensor products, \begin{equation}
{(I}_{E^{\otimes k}}\otimes\eta)(\varphi_{k}(a)\otimes I_{H})={(\varphi}_{k+1}(a)\otimes I_{H}){(I}_{E^{\otimes k}}\otimes\eta).\label{eq:balanced}\end{equation}
Consequently, we may define a map $U:\mathcal{F}(E^{\sigma})\otimes_{\iota}H\to\mathcal{F}(E)\otimes_{\sigma}H$,
where $\iota$ denotes the identity representation of $\sigma(M)'$
in $B(H)$, so that on elements of the form $\eta_{1}\otimes\eta_{2}\otimes\cdots\eta_{n}\otimes h\in\mathcal{F}(E^{\sigma})\otimes_{\iota}H$,
$U$ is given by the formula \begin{equation}
U(\eta_{1}\otimes\eta_{2}\otimes\cdots\eta_{n}\otimes h)=(I_{E^{\otimes(n-1)}}\otimes\eta_{1})(I_{E^{\otimes(n-2)}}\otimes\eta_{2})\cdots(I_{E}\otimes\eta_{n-1})\eta_{n}h\text{.}\label{eq:FourierTransform}\end{equation}
In \cite{MS05}, we called $U$ the \emph{(inverse) Fourier transform}
mapping $\mathcal{F}(E^{\sigma})\otimes_{\iota}H$ to $\mathcal{F}(E)\otimes_{\sigma}H$
determined by $E$ and $\sigma$. It plays a fundamental role in our
theory, as demonstrated by the following theorem, which is a restatement
of parts of Lemma 3.8 and Theorem 3.9 of \cite{MS04}.

\begin{theorem}\label{Commutants}Let $\sigma:M\to B(H)$ be a faithful
normal representation of $M$ on the Hilbert space $H$. Then the
inverse Fourier transform $U:\mathcal{F}(E^{\sigma})\otimes_{\iota}H\to\mathcal{F}(E)\otimes_{\sigma}H$
is a Hilbert space isomorphism such that the map\begin{equation}
X\to U\iota^{\mathcal{F}(E^{\sigma})}(X)U^{*}\label{eq:Ftrans2}\end{equation}
from $H^{\infty}(E^{\sigma})$ to $B(\mathcal{F}(E)\otimes_{\sigma}H)$
is an ultraweakly homeomorphic, completely isometric isomorphism from
$H^{\infty}(E^{\sigma})$ onto the commutant of $\sigma^{\mathcal{F}(E)}(H^{\infty}(E))$.
Likewise, the map \begin{equation}
X\to U^{*}\sigma^{\mathcal{F}(E)}(X)U\label{eq:Ftrans3}\end{equation}
is an ultraweakly continuous, completely isometric isomorphism from
$H^{\infty}(E)$ onto the commutant of $\iota^{\mathcal{F}(E^{\sigma})}(H^{\infty}(E^{\sigma}))$.
\end{theorem}

There is a formula for $U^{-1}$, but it is somewhat involved, as
may be seen from the proof of \cite[Corollary 3.10]{MS04} and one
of our goals is to circumvent it in calculations. Consequently, we
shall not develop it here.

The thrust of Proposition 5.1 of \cite{MS04} is that one can also
express $\widehat{X}(\eta^{*})$ in terms of the map defined by equation
(\ref{eq:Ftrans3}) and a ``Cauchy'' kernel expressed in terms of
$\eta$ that we define as follows. Write $\eta^{(n)}:H\rightarrow E^{\otimes n}\otimes_{\sigma}H$
for the operator given by the formula:\begin{equation}
\eta^{(n)}=(I_{E^{\otimes(n-1)}}\otimes\eta)(I_{E^{\otimes(n-2)}}\otimes\eta)\cdots(I_{E}\otimes\eta)\eta,\label{EtaToTheN}\end{equation}
 and write $\eta^{(0)}=I_{H}$. Clearly we have the recursive relation:\begin{equation}
\eta^{(n+1)}=(I_{E}\otimes\eta^{(n)})\eta=(I_{E^{\otimes n}}\otimes\eta)\eta^{(n)}\text{,}\label{recursive}\end{equation}
 which is a consequence of the formulas first proved in \cite[Lemmas 2.1 and 2.2]{MS99}.

\begin{definition}\label{CauchyKernel}The \emph{Cauchy kernel} defined
by an element $\eta\in\overline{\mathbb{D}(E^{\sigma})},$ $C(\eta),$
is the operator from $H$ to $\mathcal{F}(E)\otimes_{\sigma}H$ given
by the equation \[
C(\eta):=\left[\begin{array}{cccc}
\eta^{(0)}, & \eta^{(1)}, & \eta^{(2)}, & \cdots\end{array}\right]^{\intercal}\text{.}\]
\end{definition}

Observe that when the norm of $\eta$ is less than $1$, $C(\eta)$
is bounded with norm at most $\frac{1}{1-\Vert\eta\Vert}$. For $\eta$
of norm $1$, $C(\eta)$ may not make sense as a bounded operator
on $H$. Indeed, it is possible that $C(\eta)h$ makes sense only
when $h=0$. Observe, too, that from the definition of $\eta^{(n)}$
in (\ref{EtaToTheN}) and equation (\ref{eq:balanced}), we see immediately
that $C(\eta)$ is an element of $\mathcal{F}(E)^{\sigma}$ when $\eta\in\mathbb{D}(E^{\sigma})$.
The following proposition is a restatement of Proposition 5.1 of \cite{MS04}.
It is the starting point of our analysis.

\begin{proposition}\label{Prop5.1HAPaper}Let $\sigma:M\to B(H)$
be a normal representation of $M$ on a Hilbert space $H$ and let
$\eta\in\mathbb{D}(E)$ be given. Further, let $\rho$ be the representation
of $H^{\infty}(E)$ on $\mathcal{F}(E^{\sigma})\otimes_{\iota}H$
defined by equation (\ref{eq:Ftrans3}) and let $\iota_{H}$ be the
embedding of $H$ in $\mathcal{F}(E^{\sigma})\otimes_{\iota}H$ as
the zeroth summand. Then \[
\widehat{X}(\eta^{*})=C(\eta)^{*}U\rho(X)\iota_{H}.\]
\end{proposition}

\section{The Poisson Kernel}

We continue with the notation established above and let $E$ be a
fixed $W^{\ast}$-correspondence over a von Neumann algebra $M$ and
we let $\sigma$ be a normal representation on a Hilbert space $H$.

\begin{definition}\label{PoissonKernel}For $\eta$ in the closed
disc $\overline{\mathbb{D}(E^{\sigma})}$, we write $\Delta_{*}(\eta):=(I_{H}-\eta^{\ast}\eta)^{\frac{1}{2}}$
and we define the \emph{Poisson kernel,} $K(\eta)$, by the formula,\[
K(\eta)=(I_{\mathcal{F}(E)}\otimes\Delta_{*}(\eta))C(\eta)=(I_{\mathcal{F}(E)}\otimes\Delta_{*}(\eta))\left[\begin{array}{cccc}
\eta^{(0)}, & \eta^{(1)}, & \eta^{(2)}, & \cdots\end{array}\right]^{\intercal},\]
mapping $H$ to $\mathcal{F}(E)\otimes_{\sigma}H$. \end{definition}

\begin{remark} \label{FESigma} We note first that while $C(\eta)$
does not in general make sense as a bounded operator for $\eta$'s
with norm one, we shall see in a minute that $K(\eta)$ does. We note,
too, that $\Delta_{*}(\eta)$ commutes with $\sigma(M)$ and so $I_{\mathcal{F}(E)}\otimes\Delta_{*}(\eta)$
commutes with $\sigma^{\mathcal{F}(E)}(H^{\infty}(E))=\{ X\otimes I_{H}\mid X\in H^{\infty}(E)\}$.
Consequently, like the Cauchy kernel, $C(\eta)$, the Poisson kernel
$K(\eta)$ lies in $\mathcal{F}(E)^{\sigma}$. It will be useful to
recall that $\mathcal{F}(E)^{\sigma}$ is a $W^{*}$-correspondence
over $\sigma(M)'$. Since the action of $\sigma(M)'$ on $H$ is given
by the identity representation $\iota$, we shall denote the left
action of $\sigma(M)'$ on $\mathcal{F}(E)^{\sigma}$ by $\varphi_{\infty,\iota}$
to distinguish it from $\varphi_{\infty}$. Likewise, we write $\varphi_{\iota}$
and $\varphi_{k,\iota}$ to distinguish between the representations
induced from $\sigma$ and those induced from $\iota$. So for $c\in\sigma(M)'$
and $\eta\in\mathcal{F}(E)^{\sigma}$, $\varphi_{\infty,\iota}(c)\eta={(I}_{\mathcal{F}(E)}\otimes c)\eta$.
In particular, we may write $K(\eta)=\varphi_{\infty,\iota}(\Delta_{*}(\eta))C(\eta)$.
%
\begin{comment}
can also be written as $ $, whereObserve that $\eta^{(n)}\in(E^{\otimes n})^{\sigma}$,
so $C(\eta)$ lies in $\mathcal{F}(E)^{\sigma}$. We shall refer to
$C(\eta)$ either as a map from $H$ to $\mathcal{F}(E)\otimes_{\sigma}H$
or as an element of $\mathcal{F}(E)^{\sigma}$. The left-action map
of $\sigma(M)'$ on $\mathcal{F}(E)^{\sigma}$ will be denoted $\varphi_{\infty,\sigma}$
and, for $a\in\sigma(M)'$, $\varphi_{\infty,\sigma}(a)C(\eta)$ is
$(I_{\mathcal{F}(E)}\otimes c)\circ C(\eta)$.
\end{comment}
{} The inner product on $\mathcal{F}(E)^{\sigma}$ is simply $\langle X,Y\rangle=X^{*}Y$.
So, for $a\in\sigma(M)'$ and $\eta,\zeta\in\mathbb{D}(E^{\sigma})$,
\[
\langle C(\eta),\varphi_{\infty,\iota}(a)C(\zeta)\rangle=\sum_{k}\langle\eta^{\otimes k},\varphi_{k,\iota}(a)\zeta^{\otimes k}\rangle=\sum_{k}\theta_{\eta,\zeta}^{k}(a)=(id-\theta_{\eta,\zeta})^{-1}(a)\]
 where $\theta_{\eta,\zeta}(a)=\langle\eta,\varphi_{\iota}(a)\zeta\rangle$.
%
\begin{comment}
Therefore, if we write $k(\eta,\zeta)=(id-\theta_{\eta,\zeta})^{-1}$
and $S=\mathbb{D}(E^{\sigma})$, we see that $k(\cdot,\cdot)$ is
a completely positive definite kernel on $S$ from $\sigma(M)'$ to
itself in the sense of \cite[Definition 3.2.2]{BBLS04} and the $\sigma(M)'$-correspondence
whose existence is proved in \cite[Theorem 3.2.3]{BBLS04} is the
sub $W^{*}$-correspondence of $\mathcal{F}(E)^{\sigma}$ generated
by $\{ C(\eta):\eta\in\mathbb{D}(E)\}$. We write $F(k)$ for it and
view it as the \emph{reproducing kernel correspondence} associated
with $k$%
\marginpar{Is reference to $F(k)$ necessary?%
}.
\end{comment}
{}

\end{remark}

%
\begin{comment}
For $\eta\in E^{\sigma}$, we We next define the \emph{Cauchy kernel
determined by} Note that when we write $M\otimes_{\sigma}H$ for the
zeroth summand in $\mathcal{F}(E)\otimes_{\sigma}H$, we identify
it with $H$, as we may, via $\sigma$. Note, too, that since $\Vert\eta\Vert<1$,
$C(\eta)$ is a bounded operator from $H$ to $\mathcal{F}(E)\otimes_{\sigma}H$
of norm at most $\frac{1}{1-\Vert\eta\Vert}$.
\end{comment}
{}

\begin{proposition} \label{isometry}For all $\eta\in\overline{\mathbb{D}(E^{\sigma})}$,
$K(\eta)$ is a contraction mapping $H$ to $\mathcal{F}(E)\otimes_{\sigma}H$.
If $||\eta||<1$, then $K(\eta)$ is an isometry.\end{proposition}

\begin{proof}\begin{multline}
K(\eta)^{\ast}K(\eta)  =C(\eta)^{\ast}(I_{\mathcal{F}(E)}\otimes(\Delta(\eta))^{2})C(\eta)\\
=\left[\begin{array}{cccc}
\eta^{(0)\ast}, & \eta^{\ast}, & \eta^{(2)\ast}, & \cdots\end{array}\right](I_{\mathcal{F}(E)}\otimes(\Delta(\eta))^{2})\left[\begin{array}{cccc}
\eta^{(0)}, & \eta, & \eta^{(2)}, & \cdots\end{array}\right]^{T} \\
=\sum_{n}\eta^{(n)\ast}(I_{E^{\otimes n}}\otimes\Delta(\eta))^{2}\eta^{(n)} \\
=\lim_{N\rightarrow\infty}\sum_{n=0}^{N}\eta^{(n)\ast}\eta^{(n)}-\eta^{(n+1)\ast}\eta^{(n+1)}=I_{H}-\lim_{N\rightarrow\infty}\eta^{(N+1)\ast}\eta^{(N+1)}.\label{eq:Keta}\end{multline}
 The passage from the third line to the fourth is a consequence of
equation \ref{recursive}. Also, by \ref{recursive} and the fact
that the norm of $\eta$ is at most one, the sequence $\{\eta^{(N+1)\ast}\eta^{(N+1)}\}_{N\geq0}$
is a decreasing sequence of contractions on $H$. Therefore the sequence
${\{ I}_{H}-\eta^{(N+1)\ast}\eta^{(N+1)}\}$ converges strongly to
a contraction on $H$. The limit is $I_{H}$ if $\eta$ is a strict
contraction. \end{proof}

%
\begin{comment}
If $\eta\in\mathbb{D}(E^{\sigma})$, then $\eta^{\ast}$, which is
an operator from $E\otimes_{\sigma}H$ to $H$ determines a contractive
bimodule map, $\widehat{\eta^{\ast}}$, from $E$ to $B(H)$ via the
formula $\widehat{\eta^{\ast}}(\xi)h=\eta^{\ast}(\xi\otimes h)$.
The bimodule property means $\widehat{\eta^{\ast}}(\varphi(a)\xi b)=\sigma(a)\widehat{\eta^{\ast}}(\xi)\sigma(b)$
for all $a,b\in M$, and all $\xi\in E$.
\end{comment}
{}The following lemma shows that the values of the Poisson kernel are
``operator eigenvectors'' for the adjoints of the creation operators.
The ``operator eigenvalue'' for $T_{\xi}^{\ast}\otimes I$ determined
by $\eta\in\mathbb{D}(E^{\sigma})$ is $\widehat{\eta^{\ast}}(\xi)^{\ast}$.

\begin{lemma} \label{eigenelement}For all $\xi\in E$ and all $\eta\in\mathbb{D}(E^{\sigma})$,\[
(T_{\xi}^{\ast}\otimes I)K(\eta)=K(\eta)\widehat{\eta^{\ast}}(\xi)^{\ast}\text{.}\]
 \end{lemma}

\begin{proof} Since $\Vert\eta\Vert<1$, the operator on $H$, $\Delta_{*}(\eta)$,
is invertible. Also, $I\otimes\Delta_{*}(\eta)$ commutes with $(T_{\xi}\otimes I)^{*}$
so it suffices to prove that $\widehat{\eta^{\ast}}(\xi)C(\eta)^{\ast}=C(\eta)^{\ast}(T_{\xi}\otimes I)$
as operators from $\mathcal{F}(E)\otimes_{\sigma}H$ to $H$. To prove
equality, it suffices to evaluate both sides on an element of the
form $\zeta\otimes h\in E^{\otimes n}\otimes H$. By definition of
$C(\eta)$ and the formula (\ref{EtaToTheN}),\begin{align*}
C(\eta)^{\ast}(\zeta\otimes h) & =\eta^{(n)\ast}(\zeta\otimes h)\\
 & =\eta^{\ast}(I_{E}\otimes\eta)^{\ast}\cdots(I_{E^{\otimes(n-1)}}\otimes\eta)^{\ast}(\zeta\otimes h)\text{.}\end{align*}
 Consequently, \begin{align*}
\widehat{\eta^{\ast}}(\xi)C(\eta)^{\ast}(\zeta\otimes h) & =\widehat{\eta^{\ast}}(\xi)\eta^{\ast}(I_{E}\otimes\eta)^{\ast}\cdots(I_{E^{\otimes(n-1)}}\otimes\eta)^{\ast}(\zeta\otimes h)\\
 & =\eta^{\ast}(\xi\otimes(\eta^{\ast}(I_{E}\otimes\eta)^{\ast}\cdots(I_{E^{\otimes(n-1)}}\otimes\eta)^{\ast}(\zeta\otimes h))\\
 & =\eta^{\ast}(I_{E}\otimes\eta)^{\ast}\cdots(I_{E^{\otimes(n-1)}}\otimes\eta)^{\ast}(I_{E^{\otimes n}}\otimes\eta)^{\ast}(\xi\otimes\zeta\otimes h)\\
 & =C(\eta)^{\ast}(\xi\otimes\eta\otimes h)\\
 & =C(\eta)^{\ast}(T_{\xi}\otimes I)(\zeta\otimes h)\text{.}\end{align*}

\end{proof}

\begin{theorem} \label{Reproducing1}For all $\eta\in\mathbb{D}(E^{\sigma})$
and all $X\in H^{\infty}(E)$,\begin{equation}
K(\eta)\widehat{X}(\eta^{*})^{*}=(X^{*}\otimes I_{H})K(\eta)\label{Poisson 0}\end{equation}
 and \begin{equation}
\widehat{X}(\eta^{\ast})=K(\eta)^{\ast}(X\otimes I)K(\eta)\text{.}\label{Poisson 1}\end{equation}

\end{theorem}

\begin{proof} Remark \ref{FESigma} and Proposition \ref{eigenelement}
show that formula (\ref{Poisson 0}) holds for all $X$ of the form
$X=T_{\xi}$ and $X=\varphi_{\infty}(a)$, $\xi\in E$ and $a\in M$.
(Note that $\hat{\eta^{*}}(\xi)=\hat{T_{\xi}}(\eta^{*})$). Further
these two results show that the range of $K(\eta)$ is invariant under
all these operators. Thus the formula holds for the ultraweakly closed
algebra of operators generated by all the $T_{\xi}$ and all the $\varphi_{\infty}(a)$,
$\xi\in E$ and $a\in M$. Thus the formula (\ref{Poisson 0}) holds
for all $X\in H^{\infty}(E)$. See the discussion on page 384 of \cite{MS04}
and \cite[Corollary 2.14]{MS04}. Equation (\ref{Poisson 1}) follows
from (\ref{Poisson 0}) since $K(\eta)$ is an isometry. \end{proof}

\begin{remark}\label{Dilation}Formula \ref{Poisson 1} gives another
proof that the minimal isometric dilation of the representation of
$H^{\infty}(E)$ on $H$ determined by $\eta$ in the open disc $\mathbb{D}(E^{\sigma})$
is an induced representation of $H^{\infty}(E)$ acting on $\mathcal{F}(E)\otimes_{\sigma}H$:
$X\to X\otimes I$ \cite[Theorem 2.13]{MS04}.

\end{remark}

%
\begin{comment}
Recall the definition of the Hilbert space isomorphism $U$ from $\mathcal{F}(E^{\sigma})\otimes_{\iota}H$
to $\mathcal{F}(E)\otimes_{\sigma}H$, introduced in \cite[Lemma 3.8]{MS04}.
Theorem 3.9 of \cite{MS04} says that if we define $\rho:H(E^{\sigma})\rightarrow B(\mathcal{F}(E)\otimes_{\sigma}H)$
by the formula $\rho(X)=U(X\otimes I_{H})U^{-1}$, then $\rho$ is
an isometric, ultraweakly continuous isomorphism from $H^{\infty}(E^{\sigma})$
onto the commutant of $\sigma^{\mathcal{F}(E)}(H^{\infty}(E))$.
\end{comment}
{}The following theorem is our replacement for \cite[Proposition 5.1]{MS04}.
It captures more clearly the roles played by the various constructs.
We let $\iota_{H}$ denote the embedding of $H$ into $\mathcal{F}(E)\otimes_{\sigma}H$,
and we write $P_{H}$ for its adjoint. Also, $\rho$ is the representation
of $H^{\infty}(E^{\sigma})$ defined in equation (\ref{eq:Ftrans2}).

\begin{theorem} \label{Reproducing2}For all $\eta\in\mathbb{D}(E^{\sigma})$
and all $X\in H^{\infty}(E)$, $K(\eta)=\rho(\Delta(\eta)(I-T_{\eta})^{-1})\iota_{H}$,
and \begin{align*}
\hat{X}(\eta^{\ast}) & =K(\eta)^{\ast}(X\otimes I)K(\eta)\\
 & =P_{H}\rho(\Delta_{*}(\eta)(I-T_{\eta})^{-1}))^{\ast}(X\otimes I_{H})\rho(\Delta_{*}(\eta)(I-T_{\eta})^{-1}))\iota_{H}\\
 & =P_{H}\{\rho((I-T_{\eta})^{-1})^{\ast}(I_{\mathcal{F}(E)}\otimes\Delta_{*}(\eta)^{2})\rho((I-T_{\eta})^{-1}))\}(X\otimes I_{H})\iota_{H}\text{.}\end{align*}
 \end{theorem}

\begin{proof} Since $I_{\mathcal{F}(E)}\otimes\Delta_{*}(\eta)=\rho(\Delta_{*}(\eta))$
by \cite[Theorem 3.9]{MS04}, it suffices to prove that $C(\eta)=\rho((I-T_{\eta})^{-1})\iota_{H}$.
%
\begin{comment}
This could be deduced from \cite[Proposition 5.1]{MS04}, but the
following calculations are a bit clearer.
\end{comment}
{} Since $({I-T}_{\eta})^{-1}=\sum_{n=0}^{\infty}T_{\eta}^{n}$, it
suffices to note that for $h\in H$, $\rho(T_{\eta}^{n})h=U(\eta\otimes\eta\cdots\otimes\eta\otimes h)=(I_{E^{\otimes(n-1)}}\otimes\eta)(I_{E^{\otimes(n-2)}}\otimes\eta)\cdots(I_{E}\otimes\eta)\eta h=\eta^{(n)}h$.
\end{proof}

\vspace{4mm}

\section{Characteristic Operator Functions and Canonical Models}

In \cite{MS05} we studied canonical models for representations of
the Hardy algebras. So, given $\eta\in\mathbb{D}(E^{\sigma})$, it
makes sense and is of interest to investigate how the canonical model
of the representation $\sigma\times\widehat{\eta^{*}}$ is related
to the Poisson kernel $K(\cdot)$. We shall see that they are closely
related. We fix $\eta\in\mathbb{D}(E^{\sigma})$ for the rest of this
section and in the computations that follow, we write $\Delta_{*}=\Delta_{*}(\eta)$,
which recall is $(I_{H}-\eta^{*}\eta)^{1/2}$, and we write $\Delta=\Delta(\eta):=(I_{E\otimes H}-\eta\eta^{*})^{1/2}$
for the defect operators associated with $\widehat{\eta^{*}}$. Note
that since $\eta$ has norm strictly less than one, the operators
$\Delta$ and $\Delta_{*}$ are invertible. Therefore their ranges
are all of $E\otimes H$ and $H$, respectively. Nevertheless, to
be consistent with the literature, we continue to denote the range
of $\Delta$ by $\mathcal{D}$ and the range of $\Delta_{*}$ by $\mathcal{D}_{*}$.
We already have noted that $\Delta_{*}$ commutes with $\sigma(M)$
and it is immediate that $\Delta$ commutes with $\varphi(M)\otimes I_{H}$.
The characteristic operator of $\widehat{\eta^{*}}$ (or, of $(\widehat{\eta^{*}},\sigma)$)
is defined in \cite[Equation (12)]{MS05} to be an operator $\Theta_{\widehat{\eta^{*}}}:\mathcal{F}(E)\otimes_{\rho}\mathcal{D}\rightarrow\mathcal{F}(E)\otimes_{\rho}\mathcal{D}_{*}$
whose complete development need not be rehearsed here (in particular,
the subscript $\rho$ in the notation need not concern us). We will
give a different definition whose equivalence with the one in \cite{MS05}
will follow easily from the next lemma. It will have the advantage
that it leads immediately to a matrix representation that is useful
for our purposes. To simplify notation, we shall write $\Theta_{\eta}$
for $\Theta_{\widehat{\eta^{*}}}$.

%
\begin{comment}
We have the following.

\begin{lemma} \label{Y}(\cite[Lemma 3.12]{MS05}) The characteristic
operator $\Theta_{\eta}$ is a contraction that satisfies the equations
\begin{equation}
(\varphi_{\infty}(a)\otimes I_{\mathcal{D}_{\ast}})\Theta_{\eta}=\Theta_{\eta}(\varphi_{\infty}(a)\otimes I_{\mathcal{D}}),\;\;\; a\in M\label{Ymod}\end{equation}
 and \begin{equation}
\Theta_{\eta}(T_{\xi}\otimes I_{\mathcal{D}})=(T_{\xi}\otimes I_{\mathcal{D}_{\ast}})\Theta_{\eta},\;\;\xi\in E.\label{commute}\end{equation}
 In particular, $\Theta_{\eta}$ is determined by its restriction
to $\mathcal{D}$ ($=M\otimes\mathcal{D}\subseteq\mathcal{F}(E)\otimes_{\rho}\mathcal{D}$).
\end{lemma}
\end{comment}
{}

\begin{lemma} \label{Y} For $i=1,2$ let $\sigma_{i}$ be a faithful
normal representation of $M$ on the Hilbert space $\mathcal{E}_{i}$
and let $Y$ be a bounded linear transformation mapping $\mathcal{F}(E)\otimes_{\sigma_{1}}\mathcal{E}_{1}\rightarrow\mathcal{F}(E)\otimes_{\sigma_{2}}\mathcal{E}_{2}$.
If $Y$ intertwines $\sigma_{1}^{\mathcal{F}(E)}$ and $\sigma_{2}^{\mathcal{F}(E)}$,
then $Y$ is completely determined by its values on $\mathcal{E}_{1}$.
Conversely, given an operator $Y_{0}$ from $\mathcal{E}_{1}$ to
$\mathcal{F}(E)\otimes_{\sigma_{2}}\mathcal{E}_{2}$, the formula\[
Y(\xi\otimes e)=\xi\otimes Y_{0}e,\]
 $\xi\otimes e\in\mathcal{F}(E)\otimes_{\sigma_{1}}\mathcal{E}_{1}$
defines a bounded operator $Y:\mathcal{F}(E)\otimes_{\sigma_{1}}\mathcal{E}_{1}\rightarrow\mathcal{F}(E)\otimes_{\sigma_{2}}\mathcal{E}_{2}$
that intertwines $\sigma_{1}^{\mathcal{F}(E)}$ and $\sigma_{2}^{\mathcal{F}(E)}$.\end{lemma}

\begin{proof}The proof is immediate from Theorem \ref{Commutants}.
The only thing that might be at issue is how to handle different spaces
and different representations of $M$, $(\sigma_{i},\mathcal{E}_{i})$,
$i=1,2$. One simply forms the direct sum of $\sigma_{1}$ and $\sigma_{2}$
and induces that. Operators on the resulting space $\mathcal{F}(E)\otimes(\mathcal{E}_{1}\oplus\mathcal{E}_{2})=\mathcal{F}(E)\otimes(\mathcal{E}_{1})\oplus\mathcal{F}(E)\otimes(\mathcal{E}_{2})$
have a $2\times2$ matrix representation, and operators that intertwine
$\sigma_{1}^{\mathcal{F}(E)}$ and $\sigma_{2}^{\mathcal{F}(E)}$
can be realized as matrices of the form $\left(\begin{array}{cc}
0 & 0\\
Y & 0\end{array}\right)$.\end{proof}

To define the characteristic operator, $\Theta_{\eta}$, determined
by an element $\eta\in\mathbb{D}(E^{\sigma})$, we note that the analysis
found in \cite[pp. 429-430]{MS05} shows that the operator $\theta_{\eta}$
defined on $\mathcal{D}$ by the formula, \begin{equation}
\theta_{\eta}d=-\eta^{*}d+(I_{1}\otimes\Delta_{*})\Delta d+\sum_{k=2}^{\infty}(I_{k}\otimes\Delta_{*})(I_{1}\otimes\eta^{(k-2)})\Delta d\label{thetad}\end{equation}
 for $d\in\mathcal{D}$, is a bounded linear operator from $\mathcal{D}$
to $\mathcal{F}(E)\otimes_{\sigma}\mathcal{D}_{*}$.

\begin{definition}

For $\eta\in E^{\sigma}$, the characteristic operator determined
by $\eta$ is the operator $\Theta_{\eta}:\mathcal{F}(E)\otimes_{\varphi\otimes I\vert\mathcal{D}}\mathcal{D}\to\mathcal{F}(E)\otimes_{\sigma}\mathcal{D}_{*}$defined
by the formula\begin{equation}
\Theta_{\eta}(\xi\otimes d)=\xi\otimes\theta_{\eta}d,\label{eq:td}\end{equation}
for $d\in D$ and $\xi\in\mathcal{F}(E)$.\end{definition}

%
\begin{comment}
The characteristic operator $\Theta_{\eta}$ is a contraction that
satisfies the equations \begin{equation}
(\varphi_{\infty}(a)\otimes I_{\mathcal{D}_{\ast}})\Theta_{\eta}=\Theta_{\eta}(\varphi_{\infty}(a)\otimes I_{\mathcal{D}}),\;\;\; a\in M\label{Ymod}\end{equation}
 and \begin{equation}
\Theta_{\eta}(T_{\xi}\otimes I_{\mathcal{D}})=(T_{\xi}\otimes I_{\mathcal{D}_{\ast}})\Theta_{\eta},\;\;\xi\in E.\label{commute}\end{equation}
 In particular, $\Theta_{\eta}$ is determined by its restriction
to $\mathcal{D}$ ($=M\otimes\mathcal{D}\subseteq\mathcal{F}(E)\otimes_{\rho}\mathcal{D}$).
\end{lemma}

We write $\theta_{\eta}$ for $\Theta_{\eta}|\mathcal{D}$, the restriction
of $\Theta_{\eta}$ to $\mathcal{D}$.
\end{comment}
{}

Our next objective is to prove the following theorem which is the
principal result of this section. It was inspired in part by Popescu's
analysis in \cite{gP01} and \cite{gP06}. See \cite[Theorem 3.2]{gP06}.

\begin{theorem} \label{thetak}For $\eta\in\mathbb{D}(E^{\sigma})$,
the Poisson kernel $K(\eta)$ and the characteristic operator $\Theta_{\eta}$
are related by the equation\[
I=K(\eta)K(\eta)^{*}+\Theta_{\eta}\Theta_{\eta}^{*}\]
 on $\mathcal{F}(E)\otimes_{\sigma}\mathcal{D}_{*}$. \end{theorem}

\begin{proof}With respect to the decompositions $\mathcal{F}(E)\otimes\mathcal{D}=\mathcal{D}\oplus E\otimes\mathcal{D}\oplus E^{\otimes2}\otimes\mathcal{D}\oplus\ldots$
and $\mathcal{F}(E)\otimes\mathcal{D}_{*}=\mathcal{D}_{*}\oplus E\otimes\mathcal{D}_{*}\oplus E^{\otimes2}\otimes\mathcal{D}_{*}\oplus\ldots$,
$\Theta_{\eta}$ can be written in a matricial form $\Theta_{\eta}=(\Theta_{i,j})_{i,j=0}^{\infty}$
where $\Theta_{i,j}:E^{\otimes j}\otimes\mathcal{D}\rightarrow E^{\otimes i}\otimes\mathcal{D}_{*}$.
It follows from (\ref{eq:td}) that, for $i<j$, $\Theta_{i,j}=0$.
For $i=j$, we have $\Theta_{j,j}=I_{j}\otimes(-\eta^{*})$ and, for
$i>j$, $\Theta_{i,j}=I_{j}\otimes(I_{i-j}\otimes\Delta_{*})(I_{1}\otimes\eta^{(i-j-1)})\Delta|\mathcal{D}$.
This enables us to write the matricial form of $\Theta_{\eta}\Theta_{\eta}^{*}$
(with respect to the decomposition $\mathcal{F}(E)\otimes\mathcal{D}_{*}=\mathcal{D}_{*}\oplus E\otimes\mathcal{D}_{*}\oplus E^{\otimes2}\otimes\mathcal{D}_{*}\oplus\ldots$).
We start with the diagonal entries. \[
(\Theta_{\eta}\Theta_{\eta}^{*})_{k,k}=\sum_{l=0}^{k}\Theta_{k,l}\Theta_{k,l}^{*}=\]
 \[
I_{k}\otimes\eta^{*}\eta+\sum_{m=1}^{k}I_{k-m}\otimes(I_{m}\otimes\Delta_{*})(I_{1}\otimes\eta^{(m-1)})\Delta^{2}(I_{1}\otimes\eta^{(m-1)*})(I_{m}\otimes\Delta_{*}).\]
 But $(I_{1}\otimes\eta^{(m-1)})\Delta^{2}(I_{1}\otimes\eta^{(m-1)*})=(I_{1}\otimes\eta^{(m-1)})(I_{E\otimes H}-\eta\eta^{*})(I_{1}\otimes\eta^{(m-1)*})=I_{1}\otimes\eta^{(m-1)}\eta^{(m-1)*}-\eta^{(m)}\eta^{(m)*}$
and we get \[
(\Theta_{\eta}\Theta_{\eta}^{*})_{k,k}=I_{k}\otimes\eta^{*}\eta+\sum_{m=1}^{k}I_{k-m+1}\otimes(I_{m-1}\otimes\Delta_{*})\eta^{(m-1)}\eta^{(m-1)*}(I_{m-1}\otimes\Delta_{*})-\]
 \[
\sum_{m=1}^{k}I_{k-m}\otimes(I_{m}\otimes\Delta_{*})\eta^{(m)}\eta^{(m)*}(I_{m}\otimes\Delta_{*})=I_{k}\otimes\eta^{*}\eta+I_{k}\otimes\Delta_{*}^{2}-\]
 \[
(I_{k}\otimes\Delta_{*})\eta^{(k)}\eta^{(k)*}(I_{k}\otimes\Delta_{*})=I_{E^{\otimes k}\otimes H}-(I_{k}\otimes\Delta_{*})\eta^{(k)}\eta^{(k)*}(I_{k}\otimes\Delta_{*}).\]
 Now, fix $l<k$. Then $(\Theta_{\eta}\Theta_{\eta}^{*})_{k,l}=\sum_{m=0}^{l}\Theta_{k,m}\Theta_{l,m}^{*}$.
When $m=l$ we get \[
\Theta_{k,l}\Theta_{l,l}^{*}=I_{l}\otimes(I_{k-l}\otimes\Delta_{*})(I_{1}\otimes\eta^{(k-l-1)})\Delta(-\eta)=\]
 \[
-I_{l}\otimes(I_{k-l}\otimes\Delta_{*})(I_{1}\otimes\eta^{(k-l-1)})\eta\Delta_{*}=-I_{l}\otimes(I_{k-l}\otimes\Delta_{*})\eta^{(k-l)}\Delta_{*}.\]
 For $m<l$, $\Theta_{k,m}\Theta_{l,m}^{*}=$ $$
I_{m}\otimes(I_{k-m}\otimes\Delta_{*})(I_{1}\otimes\eta^{(k-m-1)})\Delta^{2}(I_{1}\otimes\eta^{(l-m-1)*})(I_{l-m}\otimes\Delta_{*}).$$
But $(I_{1}\otimes\eta^{(k-m-1)})\Delta^{2}(I_{1}\otimes\eta^{(l-m-1)*})=I_{1}\otimes\eta^{(k-m-1)}\eta^{(l-m-1)*}-\eta^{(k-m)}\eta^{(l-m)*}$.
Hence \[
\Theta_{k,m}\Theta_{l,m}^{*}=I_{m}\otimes(I_{k-m}\otimes\Delta_{*})(I_{1}\otimes\eta^{(k-m-1)}\eta^{(l-m-1)*})(I_{l-m}\otimes\Delta_{*})-\]
 \[
I_{m}\otimes(I_{k-m}\otimes\Delta_{*})(I_{1}\otimes\eta^{(k-m)}\eta^{(l-m-)*})(I_{l-m}\otimes\Delta_{*}).\]
 Thus \[
(\Theta_{\eta}\Theta_{\eta}^{*})_{k,l}=-I_{l}\otimes(I_{k-l}\otimes\Delta_{*})\eta^{(k-l)}\Delta_{*}+I_{l-1}\otimes(I_{k-l+1}\otimes\Delta_{*})(I_{1}\otimes\eta^{(k-l)})(I_{1}\otimes\Delta_{*})-\]
 \[
(I_{k}\otimes\Delta_{*})\eta^{(k)}\eta^{(l)*}(I_{l}\otimes\Delta_{*})=-(I_{k}\otimes\Delta_{*})\eta^{(k)}\eta^{(l)*}(I_{l}\otimes\Delta_{*}).\]
 It is easy to check, using the definition of $K(\eta)$, that the
matricial form of $K(\eta)K(\eta)^{*}$ is \[
(K(\eta)K(\eta)^{*})_{k,l}=(I_{k}\otimes\Delta_{*})\eta^{(k)}\eta^{(l)*}(I_{l}\otimes\Delta_{*})\]
 and we conclude

 \[
\Theta_{\eta}\Theta_{\eta}^{*}+K(\eta)K(\eta)^{*}=I.\]
 \end{proof}

\section{Point Evaluations on $\mathbb{D}(E)$}

Recall from \cite[Theorem 3.6]{MS04} that there is a natural isomorphism
between $E$ and ${(E}^{\sigma})^{\iota}$, where $\iota$ denotes
the identity representation of $\sigma(M)'$ on $H$. Thus we may
identify $E$ and ${(E}^{\sigma})^{\iota}$ and view elements of $H^{\infty}(E^{\sigma})$
as functions on $\mathbb{D}(E)$. This will help to shed some light
on the relation between our work and \cite{ADD90} and it will enable
us to (anti)represent $H^{\infty}(E)$ in the algebra of completely
bounded maps on $M$, $CB(M).$ For this purpose, we adopt the convention
that when $X\in H^{\infty}(E)$ and when we write $X1$, $1$ is understood
to be the identity of $M$ viewed as a vector of $\mathcal{F}(E)=M\oplus E\oplus\ldots$.
So $X1\in\mathcal{F}(E)$. We write $C(\xi)$ and $K(\xi)$, for $\xi\in\mathbb{D}(E)$,
using the obvious modifications of Definitions \ref{CauchyKernel}
and \ref{PoissonKernel}, and note that $K(\xi)=\varphi_{\infty}(\Delta_{*})C(\xi)$
where $\Delta_{*}=(I-\langle\xi,\xi\rangle)^{1/2}$. Also, we write
$\mathbb{E}_{0}$ for the conditional expectation of $H^{\infty}(E)$.
This map is defined as $\Phi_{0}$ on page 336 of \cite{MS04}. It
picks off the zeroth coefficient of an element $X\in H^{\infty}(E)$
calculated with respect to the gauge automorphism group.

\begin{theorem}\label{Phi} For $\xi\in\mathbb{D}(E)$, and $X\in H^{\infty}(E)$,
we define the map $\Phi_{X}^{\xi}:M\rightarrow M$ by the formula
\begin{equation}
\Phi_{X}^{\xi}(a)=\langle C(\xi),\varphi_{\infty}(a)X1\rangle,\label{Phi1}\end{equation}
for all $a\in M$. Then

\begin{enumerate}
\item [(1)] For each $a\in M$, $\Phi_{X}^{\xi}(a)$ is the unique element
of $M$ such that \[
(I-T_{\xi}^{*})^{-1}(\varphi_{\infty}(a)X-\varphi_{\infty}(\Phi_{X}^{\xi}(a)))\in H_{0}^{\infty}(E),\]
 where $H_{0}^{\infty}(E):={\bigvee\{ T}_{\xi}X\mid\xi\in E,{\: X\in H}^{\infty}(E)\}=H^{\infty}(E)\cap Ker(\mathbb{E}_{0})$.
\item [(2)] We have \[
X^{*}\varphi_{\infty}(a)C(\xi)=\varphi_{\infty}(\Phi_{X}^{\xi}(a)^{*})C(\xi)\]
 and, in particular, \[
X^{*}K(\xi)=\varphi_{\infty}(\Phi_{X}^{\xi}(\Delta_{*})^{*}\Delta_{*}^{-1})K(\xi).\]
 So that $K(\xi)$ is an eigenvector of $X^{*}$ (cf. Corollary~\ref{Reproducing1}).
\item [(3)] For each $\xi\in\mathbb{D}(E)$, the map $\Phi:X\mapsto\Phi_{X}^{\xi}$
is an algebra antihomomorphism from $H^{\infty}(E)$ into $CB(M)$.
\end{enumerate}
\end{theorem} \begin{proof} First note that, since $\Vert\xi\Vert<1$,
$I-T_{\xi}^{*}$ is an invertible operator on $\mathcal{F}(E)$ with
inverse equal to $I+T_{\xi}^{*}+T_{\xi}^{*2}+\ldots$ We claim that
for $X\in H^{\infty}(E)$, $(I+T_{\xi}^{*}+T_{\xi}^{*2}+\ldots)(\varphi_{\infty}(a)X-\varphi_{\infty}(\Phi_{X}^{\xi}(a)))$
lies in $H_{0}^{\infty}(E)$. If $X=T_{g}$ for some $g\in E^{\otimes n}$,
then \[
(I+T_{\xi}^{*}+T_{\xi}^{*2}+\ldots)(\varphi_{\infty}(a)X-\varphi_{\infty}(\Phi_{X}^{\xi}(a)))=\]
 \[
(I+T_{\xi}^{*}+T_{\xi}^{*2}+\ldots)(\varphi_{\infty}(a)T_{g}-\varphi_{\infty}(\langle\xi^{\otimes n},\varphi_{n}(a)g\rangle)).\]
 Note, too, that $T_{\xi}^{*k}\varphi_{\infty}(a)T_{g}=T_{\xi}^{*(k-n)}\varphi_{\infty}(\langle\xi^{\otimes n},\varphi_{n}(a)g\rangle)$,
for $k\geq n$. Thus \begin{multline}
(I+T_{\xi}^{*}+T_{\xi}^{*2}+\ldots)(\varphi_{\infty}(a)X-\varphi_{\infty}(\Phi_{X}^{\xi}(a)))\\=\varphi_{\infty}(a)T_{g}+T_{\xi}^{*}\varphi_{\infty}(a)T_{g}+T_{\xi}^{*2}\varphi_{\infty}(a)T_{g}+
\ldots+T_{\xi}^{*(n-1)}\varphi_{\infty}(a)T_{g}\in H_{0}^{\infty}(E).\end{multline}
 It follows that the result holds for all operators in a ultraweakly-dense
subalgebra of $H^{\infty}(E)$. Since the map taking $X\in H^{\infty}(E)$
to $\Phi_{X}^{\xi}(a)$ is ultraweakly-continuous, $(I+T_{\xi}^{*}+T_{\xi}^{*2}+\ldots)(\varphi_{\infty}(a)X-\varphi_{\infty}(\Phi_{X}^{\xi}(a)))$
lies in $H_{0}^{\infty}(E)$ for all $X\in H^{\infty}(E)$. To prove
uniqueness we need to show that, if $c\in\varphi_{\infty}(M)$ satisfies
$(I-T_{\xi}^{*})^{-1}c\in H^{\infty}(E)_{0}$, then $c=0$. But, since
$(I-T_{\xi}^{*})^{-1}c=(I+T_{\xi}^{*}+T_{\xi}^{*2}+\ldots)c$, this
is clear and (1) follows.

To prove (2), fix $X\in H^{\infty}(E)$, $a\in M$ and write $Y$
for $(I-T_{\xi}^{*})^{-1}(\varphi_{\infty}(a)X-\varphi_{\infty}(\Phi_{X}^{\xi}(a)))$
(in $H_{0}^{\infty}(E)$). Then $(\varphi_{\infty}(a)X-\varphi_{\infty}(\Phi_{X}^{\xi}(a)))^{*}=Y^{*}(I-T_{\xi})$.
Since $(I-T_{\xi})C(\xi)=1\in\mathcal{F}(E)$, and $Y\in H_{0}^{\infty}(E)$,
we have $(\varphi_{\infty}(a)X-\varphi_{\infty}(\Phi_{X}^{\xi}(a)))^{*}C(\xi)=Y^{*}(I-T_{\xi})C(\xi)=0$.
This, together with the observation that $K(\xi)=\varphi_{\infty}(\Delta_{*})C(\xi)$,
completes the proof of (2).

Finally, note that the linearity of the map $X\to\Phi_{X}^{\xi}$
is obvious and antimultiplicativity follows from the computation $\Phi_{XZ}^{\xi}(a)=\langle X^{*}\varphi_{\infty}(a^{*})C(\xi),Z1\rangle=\langle\varphi_{\infty}(\Phi_{X}^{\xi}(a)^{*})C(\xi),Z1\rangle=\Phi_{Z}^{\xi}(\Phi_{X}^{\xi}(a))$.
\end{proof}

\begin{remark}\label{remtheorem}

\begin{enumerate}
\item [(i)] When we fix $X$ and $\xi$ and let $a=I\in M$, we find that
$\Phi_{X}^{\xi}(I)$ is very closely related to the concept of {}``left
point evaluation\char`\"{} of $X$ at $\xi$ that was defined for
the special case of upper triangular operators in \cite{ADD90} and
studied there and in subsequent papers by various authors. (Compare
\cite[Theorem 3.3]{ADD90} with Theorem~\ref{Phi}(1)). If one adopts
the {}``reproducing kernel correspondence\char`\"{} point of view
discussed in Remark~\ref{FESigma}, this indeed can be viewed as
a point evaluation. Note, however, that the map $X\mapsto\Phi_{X}^{\xi}(I)$
is not multiplicative in general. (See also \cite[Example 2.25]{MS05a}).
%
\begin{comment}
[(ii)] In our presentation of Theorem \ref{Phi}, we fixed $\xi$
and suppressed it in the notation. If we wish to make the dependence
on $\xi$ explicit, we shall write $\Phi_{X,\xi}$ (instead of $\Phi_{X}$
) and $\Phi_{\xi}$ (for the map $X\mapsto\Phi_{X,\xi}$, instead
of $\Phi$).
\end{comment}
{}
\item [(ii)] It follows from Theorem~\ref{Phi} that, for each $\xi\in\mathbb{D}(E)$,
the kernel of the map $X\to\Phi_{X}^{\xi}$ is a two-sided ideal in
$H^{\infty}(E)$ .
\end{enumerate}
\end{remark}

\section{Curvature}

In this section we express the curvature invariant that we attached
to completely positive maps on semifinite factors \cite{MS03} in
terms of the Poisson kernel. This provides a further connection between
that work and the analysis by Popescu in \cite{gP01} and the study by Arveson \cite{wA00}. We suppose
from now on that $M$ is a semifinite factor and we fix a faithful
normal semifinite trace $\tau$ on $M$. We recall that once $\tau$
is fixed, we may define a dimension for any representation and we
can assign a natural trace to the commutant of the representation
(cf. \cite[Definition 2.1]{MS03}). Specifically, if $\sigma$ is
a normal representation of $M$ on $H$, then there is a Hilbert space
isometry $u$ from $H$ to $L^{2}(M,\tau)\otimes\ell^{2}(\mathbb{N})$,
where $L^{2}(M,\tau)$ is the $L^{2}$-space canonically associated
with $\tau$, i.e., the GNS-space, such that $u\sigma(a)=\lambda(a)\otimes I_{\ell^{2}(\mathbb{N})}u$,
for all $a\in M$, where $\lambda$ is the left representation of
$M$ on $L^{2}(M,\tau)$. The range projection of $u$, $e$, lies
in the commutant of $\lambda(M)\otimes I_{\ell^{2}(\mathbb{N})}$,
which is $\rho(M)\otimes B(\ell^{2}(\mathbb{N}))$, where $\rho$
is the right (anti) representation of $M$ on $L^{2}(M,\tau)$. The
usual trace on $\rho(M)\otimes B(\ell^{2}(\mathbb{N}))$ is $\tau\otimes tr$,
where $tr$ is the standard trace on $B(\ell^{2}(\mathbb{N}))$, i.e.,
the one that assigns to each projection in $B(\ell^{2}(\mathbb{N}))$
its rank. Then, while $u$ and $e$ are not unique, the Murray-von
Neumann equivalence class of $e$ in $\rho(M)\otimes B(\ell^{2}(\mathbb{N}))$
is uniquely determined by $\sigma$ and so, therefore, $\tau\otimes tr(e)\in[0,\infty]$
is unique. This number is called \emph{the dimension of} $H$ (or
of $\sigma$) as a module over $M$. We write this number $\dim_{\sigma}H$.
It will be important to remember, too, that the commutant of $\sigma(M)$
is spatially isomorphic to $e(\rho(M)\otimes B(\ell^{2}(\mathbb{N})))e$
via $u$ and so we can refer to the natural trace on $\sigma(M)'$
as the restriction of $\tau\otimes tr$ to $e(\rho(M)\otimes B(\ell^{2}(\mathbb{N})))e$.
We shall do this and we shall denote it by $\mathrm{tr}{}_{\sigma(M)'}$.
If $E$ is a $W^{*}$-correspondence over $M$ of the kind we have
been studying, then the \emph{(left) dimension} of $E$ is defined
to be the dimension of the representation $\varphi\otimes I$, representing
$M$ on $E\otimes_{\lambda}L^{2}(M,\tau)$. We denote this dimension
by $\dim_{l}(E)$. (See \cite[Definition 2.5]{MS03}.)

An $\eta\in\overline{\mathbb{D}(E^{\sigma})}$ defines a completely
positive map $P=P_{\eta}$ on $\sigma(M)'$ via the formula $P(a)=\eta^{*}(I_{E}\otimes a)\eta$,
$a\in\sigma(M)'$. Alternatively, given the formula for the inner
product in $E^{\sigma}$, $P(a)=\langle\eta,a\eta\rangle$. And conversely,
given a completely positive map $P$ on $\sigma(M)'$ there a $W^{*}$-correspondence
$E$ over $M$ and an $\eta\in\overline{\mathbb{D}(E^{\sigma})}$
such that $P=P_{\eta}$\cite[Corollary 2.23]{MS02}.

\begin{definition}\label{kappa}Let $E$ be a $W^{*}$-correspondence
over the von Neumann algebra $M$ with $\dim_{l}(E):=d$, and let
$\sigma$ be a representation of $M$ on the Hilbert space $H$. Then
for $\eta\in\overline{\mathbb{D}(E^{\sigma})}$, the \emph{curvature}
of $\eta$ is defined to be the curvature of $P_{\eta}$ in the sense
of \cite[Definition 3.1]{MS03}, which is the limit \[
\lim_{N\to\infty}\frac{\mathrm{tr}{}_{\sigma(M)'}(I-P_{\eta}^{N+1}(I_{H}))}{\sum_{k=0}^{N}d^{k}},\]
 and will be denoted $\kappa(\eta)$. \end{definition}

The limit exists, as was shown in \cite[Theorem 3.3]{MS03}, where
alternate formulas for $\kappa(\eta)$ may also be found. The basis
for the calculations we make here is the following lemma, whose proof
may be assembled easily from \cite{MS03}.

\begin{lemma}\label{FundCalc}Let $E$ and $F$ be $W^{*}$-correspondences
over $M$ and let $\sigma$ be normal representation of $M$ on a
Hilbert space $H$. Then

\begin{enumerate}
\item $\dim_{l}(E\otimes F)=\dim_{l}E\times\dim_{l}F$.
\item If $\eta,\zeta\in E^{\sigma}$, then \begin{equation}
\mathrm{tr}_{\sigma(M)'}{(\zeta}^{*}\eta)=\mathrm{tr}_{(\varphi(M)\otimes I_{H})'}(\eta\zeta^{*})\label{eq:trace}\end{equation}
where $\varphi$ denotes the left action of $M$ on $E$.
\item For all positive $x$ in $\sigma(M)'$, $\textrm{tr}_{(\varphi(M)\otimes I_{H})'}(I_{E}\otimes x)=\textrm{tr}_{\sigma(M)'}(x)\cdot\dim_{l}E$.
\item $\dim_{\varphi\otimes I_{H}}(E\otimes_{\sigma}H)=\dim_{l}(E)\cdot\dim_{\sigma}H$.
\end{enumerate}
\end{lemma}

\begin{proof}The first assertion is proved as Corollary 2.8 in \cite{MS03}.
The second assertion is embedded in the proof of \cite[Proposition 2.12]{MS03}.
For the sake of clarity we repeat the salient part of it here. Form
the direct sum $H\oplus(E\otimes_{\sigma}H)$ and let $\tilde{\sigma}=\sigma\oplus(\varphi\otimes I_{H})$
be the representation of $M$ acting on this space. Then the commutant
of $\tilde{\sigma}(M)$ is the set of all matrices of the form $\left(\begin{array}{cc}
a & b\\
c & d\end{array}\right)$ where $a\in\sigma(M)'$, $d\in(\varphi(M)\otimes I_{H})'$, $b\varphi(x)\otimes I_{H}=\sigma(x)b$
and $c\sigma(x)=\varphi(x)\otimes I_{H}c$ for all $x\in M$. Further,
it is easy to see that $\mathrm{tr}_{\tilde{\sigma}(M)'}\left(\begin{array}{cc}
a & b\\
c & d\end{array}\right)=\mathrm{tr}_{\sigma(M)'}(a)+\mathrm{tr}_{(\varphi(M)\otimes I_{H})'}(d)$. Thus we find that \begin{multline*} \mathrm{tr}_{\sigma(M)'}{(\zeta}^{*}\eta)=\mathrm{tr}_{\tilde{\sigma}(M)'}\left(\begin{array}{cc}
\zeta^{*}\eta & 0\\
0 & 0\end{array}\right)\\=\mathrm{tr}_{\tilde{\sigma}(M)'}\left(\left(\begin{array}{cc}
0 & \zeta^{*}\\
0 & 0\end{array}\right)\left(\begin{array}{cc}
0 & 0\\
\eta & 0\end{array}\right)\right)=\mathrm{tr}_{\tilde{\sigma}(M)'}\left(\left(\begin{array}{cc}
0 & 0\\
\eta & 0\end{array}\right)\left(\begin{array}{cc}
0 & \zeta^{*}\\
0 & 0\end{array}\right)\right)\\=\mathrm{tr}_{\tilde{\sigma}(M)'}\left(\begin{array}{cc}
0 & 0\\
0 & {\eta\zeta}^{*}\end{array}\right)=\mathrm{tr}_{(\varphi(M)\otimes I_{H})'}(\eta\zeta^{*}).\end{multline*} The third assertion is \cite[Lemma 2.7]{MS03}, and the last assertion
follows from the third by taking $x=I_{H}$.\end{proof}

Let $E$ and $\sigma$ be fixed, now, write $d$ for $\dim_{l}E$,
and write $P_{m}$ for the projection of $\mathcal{F}(E)$ onto $E^{\otimes m}$.
Also, write $P_{\leq m}$ for the sum $\sum_{k\leq m}P_{k}$. Then
it is evident from Lemma \ref{FundCalc} that \[
\textrm{tr}_{(\varphi_{\infty}(M)\otimes I_{H})'}P_{m}\otimes I_{H}=\dim_{(\varphi_{\infty}\otimes I_{H})}E^{\otimes m}\otimes_{\sigma}H=d^{m}\dim_{\sigma}H.\]

\begin{theorem}\label{CurvaturePoisson}If $\eta\in\overline{\mathbb{D}(E^{\sigma})}$,
then:
\begin{enumerate}

\item If $d:=\dim_l E$ is finite, then \begin{multline}
\kappa(\eta)  =  \lim_{N\to\infty}\frac{\sum_{k=0}^{N}\textrm{tr}_{\sigma(M)'}[K(\eta)^{*}(P_{k}\otimes I_{H})K(\eta)]}{(1+d+d^{2}+\cdots+d^{N})} \\
=  \lim_{N\to\infty}\frac{\textrm{tr}_{(\varphi_{\infty}(M)\otimes I_{H})'}[(P_{\leq N}\otimes I_{H})(K(\eta)K(\eta)^{*})(P_{\leq N}\otimes I_{H})]}{(1+d+d^{2}+\cdots+d^{N})}.\label{eq:Curv1}\end{multline}

\item If $d\geq1$, then \begin{multline}
\kappa(\eta)  =  \lim_{N\to\infty}\frac{\textrm{tr}_{\sigma(M)'}[K(\eta)^{*}(P_{N}\otimes I_{H})K(\eta)]}{d^{N}} \\
=  \lim_{N\to\infty}\frac{\textrm{tr}_{(\varphi_{\infty}(M)\otimes I_{H})'}[(P_{N}\otimes I_{H})K(\eta)K(\eta)^{*}(P_{N}\otimes I_{H})]}{d^{N}}.\label{eq:Curv2}\end{multline}

\item If $d<1$, and if $\dim_{\sigma}H$ is finite or more generally if $\textrm{tr}_{\sigma(M)'}(I_{H}-P_{\eta}(I_{H}))$ is finite,
then $\textrm{tr}_{(\varphi_{\infty}(M)\otimes I_{H})'}(K(\eta)K(\eta)^{*})$
is finite and \begin{multline}
\kappa(\eta)  =  (1-d)\cdot\textrm{tr}_{(\varphi_{\infty}(M)\otimes I_{H})'}(K(\eta)K(\eta)^{*})\\ = (1-d)\cdot\textrm{tr}_{\sigma(M)'}(K(\eta)^{*}K(\eta)).\label{eq:Curv3}\end{multline}
In particular, if $\Vert\eta\Vert<1$ and $d<1$, then $\kappa(\eta)=(1-d)\cdot\dim_{\sigma}(H)$
is independent of $\eta$.
\end{enumerate}

\end{theorem}

\begin{proof} By definition, $$\kappa(\eta)=\lim_{N\to\infty}\frac{\mathrm{tr}{}_{\sigma(M)'}(I-P_{\eta}^{N+1}(I_{H}))}{\sum_{k=0}^{N}d^{k}},$$
and by definition of $P_{\eta}$ and equation (\ref{eq:Keta}), the
numerator in the definition of $\kappa(\eta)$ is $\sum_{k=0}^{N}\textrm{tr}_{\sigma(M)'}[K(\eta)^{*}(P_{k}\otimes I_{H})K(\eta)]$.
This proves the equality of the first two terms in equation (\ref{eq:Curv1}).
The equality of the third term with the first two is immediate from
equation (\ref{eq:trace}) in Lemma \ref{FundCalc} (when $\mathcal{F}(E)$
is used in place of $E$). For the second equation, write the sum
$\sum_{k=0}^{N}\textrm{tr}_{\sigma(M)'}[K(\eta)^{*}(P_{k}\otimes I_{H})K(\eta)]$
as $\textrm{tr}_{\sigma(M)'}[K(\eta)^{*}(P_{\leq N}\otimes I_{H})K(\eta)]$,
then the two numerators in equation (\ref{eq:Curv2}) are the same
by Lemma \ref{FundCalc}. But $\lim_{N\to\infty}\frac{\sum_{k=0}^{N}\textrm{tr}_{\sigma(M)'}[K(\eta)^{*}(P_{k}\otimes I_{H})K(\eta)]}{(1+d+d^{2}+\cdots+d^{N})}=\lim_{N\to\infty}\frac{\textrm{tr}_{\sigma(M)'}[K(\eta)^{*}(P_{N}\otimes I_{H})K(\eta)]}{d^{N}}$
and \begin{eqnarray*}
\lim_{N\to\infty}\frac{\textrm{tr}_{(\varphi_{\infty}(M)\otimes I_{H})'}[(P_{\leq N}\otimes I_{H})(K(\eta)K(\eta)^{*})(P_{\leq N}\otimes I_{H})]}{(1+d+d^{2}+\cdots+d^{N})}\\
=\lim_{N\to\infty}\frac{\textrm{tr}_{(\varphi_{\infty}(M)\otimes I_{H})'}[(P_{N}\otimes I_{H})K(\eta)K(\eta)^{*}(P_{N}\otimes I_{H})]}{d^{N}}\end{eqnarray*}
using \cite[Lemma 3.2]{MS03} (first noted in \cite[Page 280]{gP01})
and the arguments from the proof of \cite[Theorem 3.3]{MS03}. This
proves equation (\ref{eq:Curv2}). Finally, for equation (\ref{eq:Curv3}),
observe that when $d<1$ and $\textrm{tr}_{\sigma(M)'}(I_{H}-P_{\eta}(I_{H}))<\infty$,
the argument in the last paragraph of the proof of \cite[Theorem 3.3]{MS03}
shows that the traces $\textrm{tr}_{\sigma(M)'}(I_{H}-P_{\eta}^{N}(I_{H}))$
increase to a finite limit. Since each of these traces equals $\textrm{tr}_{(\varphi_{\infty}(M)\otimes I_{H})'}[(P_{\leq N-1}\otimes I_{H})(K(\eta)K(\eta)^{*})(P_{\leq N-1}\otimes I_{H})]$
by Lemma \ref{FundCalc}, the normality of the trace, $\textrm{tr}_{(\varphi_{\infty}(M)\otimes I_{H})'}$,
implies that $\textrm{tr}_{(\varphi_{\infty}(M)\otimes I_{H})'}(K(\eta)K(\eta)^{*})<\infty$.
As in the proof of \cite[Theorem 3.3]{MS03}, the proof of equation
(\ref{eq:Curv3}) is immediate from the definition of $\kappa(\eta)$,
the formula for the partial sums of a geometric series, and the fact
that $d<1$. \end{proof}

Our final goal is to relate the curvature, $\kappa(\eta)$, with the
trace of the ``curvature operator'' naturally associated to $\eta$.
To define this operator, we need to say a bit more about the induced
covariant representations of $E$, Definition \ref{InducedCovariant}.
Recall that it is $(V,\varphi_{\infty}\otimes I_{H})$, where $V:E\rightarrow B(\mathcal{F}(E)\otimes_{\sigma}H)$
is defined by the formula $V(\xi)=T_{\xi}\otimes I_{H}$. The associated
map $\tilde{V}:E\otimes\mathcal{F}(E)\otimes_{\sigma}H\to\mathcal{F}(E)\otimes_{\sigma}H$
is ``simply multiplication'': $\tilde{V}(\xi\otimes(\eta\otimes h))=(\xi\otimes\eta)\otimes h$.
As we remarked earlier, while this map looks like the identity embedding
of $\sum_{k=1}^{\infty}E^{\otimes k}\otimes H$ into $\mathcal{F}(E)\otimes_{\sigma}H$,
there is, in fact, a shift involved. Specifically, if $P_{k}$ is
the projection of $\mathcal{F}(E)$ onto the summand $E^{\otimes k}$,
then a simple calculation shows that \begin{equation}
\tilde{V}(I_{E}\otimes{(P}_{k}\otimes I_{H}))={(P}_{k+1}\otimes I_{H})\tilde{V}\label{eq:Shift}\end{equation}
(see \cite[Corollary 2.4]{MS99}.) Alternatively, we may say that
$\tilde{V}^{*}$is a coisometric map in $E^{\varphi_{\infty}\otimes I_{H}}$.
We shall write $\tilde{V_{0}}:=I_{\mathcal{F}(E)\otimes_{\sigma}H}$
and recursively define $\tilde{V}_{k+1}:=\tilde{V}(I_{E}\otimes\tilde{V_{k})}$.
The map $\tilde{V}$ induces a non-unital endomorphism of ${(\varphi}_{\infty}(M)\otimes I_{H})'$
by the formula $\Phi_{V}(X)=\tilde{V}(I_{E}\otimes X)\tilde{V}^{*}$
and the powers of $\Phi_{V}$ are given by the formula $\Phi_{V}^{n}(X)=\tilde{V}_{n}(I_{E^{\otimes n}}\otimes X)\tilde{V}_{n}^{*}$\cite[Lemma 2.3]{MS99}.%
\footnote{It may be helpful to keep in mind that expressions like $I_{E}\otimes_{\sigma}X$
need not represent bounded operators unless $X$ is a (bounded) operator
in the commutant of $\sigma(M)$. That is why the formula for $\Phi_{V}$
does not make sense unless the argument is from $(\varphi_{\infty}(M)\otimes I_{H})'$.%
} We also define $\delta_{V}:{(\varphi}_{\infty}(M)\otimes I_{H})'\to{(\varphi}_{\infty}(M)\otimes I_{H})'$
by the formula $\delta_{V}(X):=X-\Phi_{V}(X)$ and we define $N:=\sum_{k=0}^{\infty}d^{-k}(P_{k}\otimes I_{H})$,
where, recall, $d=\dim_{l}E$. Thus $\Phi_{V}$, $\delta_{V}$, and
$N$ are analogues of Popescu's operators, $\phi_{S\otimes I}$, $d_{S\otimes I}$
and $N$, defined on pages 271 and 272 of \cite{gP01}. Note that
$N$ is bounded only when $d\geq1$.

\begin{definition}\label{CurvatureOperator}For $\eta\in\overline{\mathbb{D}(E^{\sigma})}$,
the \emph{curvature operator} determined by $\eta$ is defined to
be \[
\delta_{V}[K(\eta)K(\eta)^{*}]N.\]

\end{definition}

Our goal is to prove the following analogue of \cite[Theorem 2.3]{gP01}.

\begin{theorem}\label{CurvatureFormula}If $d:=\dim_{l}E\geq1$,
then for $\eta\in\overline{\mathbb{D}(E^{\sigma})}$ , \[
\kappa(\eta)=\textrm{tr}_{(\varphi_{\infty}(M)\otimes I_{H})'}\{\delta_{V}[K(\eta)K(\eta)^{*}]N\}.\]

\end{theorem}

\begin{proof}We begin by proving an analogue of \cite[Theorem 1.1]{gP01}.
For $Y\in{(\varphi}_{\infty}(M)\otimes I_{H})'$, $\sum_{k=0}^{m}\Phi_{V}^{k}(\delta_{V}(Y))=\sum_{k=0}^{m}\tilde{V_{k}}(I_{E^{\otimes k}}\otimes Y)\tilde{V}_{k}^{*}-\tilde{V}_{k+1}(I_{E^{\otimes k+1}}\otimes Y)\tilde{V}_{k+1}^{*}=Y-\tilde{V}_{m+1}(I_{E^{\otimes m+1}}\otimes Y)\tilde{V}_{m+1}^{*}$
. Since $(V,\varphi_{\infty}\otimes I_{H})$ is an induced representation
in the sense of \cite[Page 854]{MS99}, \cite[Corollary 2.10]{MS99}
implies that the ultra-strong limit, $\lim_{n\to\infty}\tilde{V}_{n}(I_{E^{\otimes n}}\otimes Y)\tilde{V}_{n}^{*}=0$.
Thus $Y=\sum_{k=0}^{\infty}\Phi_{V}^{k}(\delta_{V}(Y))$, where the
convergence is in the ultra-strong topology. Thus for each $m\geq0$,
we have on the basis of equation (\ref{eq:Shift}),\begin{multline*}
{(P}_{m}\otimes I_{H})Y(P_{m}\otimes I_{H})  =  \sum_{k\geq0}(P_{m}\otimes I_{H})\tilde{V_{k}}(I_{E^{\otimes k}}\otimes\delta_{V}(Y))\tilde{V}_{k}^{*}{(P}_{m}\otimes I_{H})\\
=  \sum_{k\geq0}\tilde{V_{k}}({I_{E^{\otimes k}}\otimes P}_{m-k}\otimes I_{H})(I_{E^{\otimes k}}\otimes\delta_{V}(Y))({I_{E^{\otimes k}}\otimes P}_{m-k}\otimes I_{H})\tilde{V}_{k}^{*}\\
=  \sum_{k=0}^{m}\tilde{V_{k}}(I_{E^{\otimes k}}\otimes({P_{m-k}\otimes I_{H})\delta}_{V}(Y)(P_{m-k}\otimes I_{H}))\tilde{V}_{k}^{*}.\end{multline*}
Thus, since $\textrm{tr}_{{(\varphi}_{\infty}(M)\otimes I_{H})'}$restricts
to $\textrm{tr}_{{(\varphi}_{m}(M)\otimes I_{H})'}$ on ${(P}_{m}\otimes I_{H})(\varphi_{\infty}(M)\otimes I_{H})'(P_{m}\otimes I_{H})$,
we see that for any operator $Y$ that has finite trace calculated with respect to $\textrm{tr}_{{(\varphi}_{\infty}(M)\otimes I_{H})'}$
and for any positive operator $Y$ in ${(\varphi}_{\infty}(M)\otimes I_{H})'$,
\begin{multline*}
\textrm{tr}_{{(\varphi}_{\infty}(M)\otimes I_{H})'}((P_{m}\otimes I_{H})Y(P_{m}\otimes I_{H}))\\  =  \sum_{k=0}^{m}\textrm{tr}_{{(\varphi}_{m}(M)\otimes I_{H})'}(\tilde{V_{k}}(I_{E^{\otimes k}}\otimes({P_{m-k}\otimes I_{H})\delta}_{V}(Y)(P_{m-k}\otimes I_{H})\tilde{V}_{k}^{*})\\
=  \sum_{k=0}^{m}\textrm{tr}_{{(\varphi}_{k}(M)\otimes I_{E^{\otimes(m-k)}\otimes H})'}(I_{E^{\otimes k}}\otimes({P_{m-k}\otimes I_{H})\delta}_{V}(Y)(P_{m-k}\otimes I_{H}))\\
=  \sum_{k=0}^{m}\textrm{tr}_{{(\varphi}_{m-k}(M)\otimes I_{H})'}(({P_{m-k}\otimes I_{H})\delta}_{V}(Y)(P_{m-k}\otimes I_{H}))d^{k}\end{multline*}
We can pass from the first line in this equation to the second and
eliminate the $\tilde{V}_{k}$, since they simply identify $E^{\otimes k}\otimes(E^{\otimes(m-k)}\otimes H)$
with $E^{\otimes m}\otimes H$, and in so doing transform the trace
on ${(\varphi}_{m}(M)\otimes I_{H})'$, $\textrm{tr}_{{(\varphi}_{m}(M)\otimes I_{H})'}$,
to the trace on ${(\varphi}_{k}(M)\otimes I_{E^{\otimes(m-k)}\otimes H})'$,
$\textrm{tr}_{{(\varphi}_{k}(M)\otimes I_{E^{\otimes(m-k)}\otimes H})'}$.
The passage to the last line is justified by part 3. of Lemma \ref{FundCalc}.
Here, $\varphi_{k}$ plays the role of $\varphi$ in the lemma, while
$\varphi_{m-k}\otimes I_{H}$ plays the role of $\sigma$. Also, of
course, part 1. of that lemma guarantees that $\dim_{l}E^{\otimes k}=d^{k}$.
So, if we divide the equation by $d^{m}$ and then change variables
in the last sum, $m-k\to k$, we conclude that \begin{multline}
\frac{\textrm{tr}_{{(\varphi}_{\infty}(M)\otimes I_{H})'}({(P}_{m}\otimes I_{H})Y(P_{m}\otimes I_{H}))}{d^{m}} \\ =  \sum_{k=0}^{m}\textrm{tr}_{{(\varphi}_{k}(M)\otimes I_{H})'}({{(P}_{k}\otimes I_{H})\delta}_{V}(Y)(P_{k}\otimes I_{H}))d^{-k} \\
=  \sum_{k=0}^{m}\textrm{tr}_{{(\varphi}_{\infty}(M)\otimes I_{H})'}(\delta_{V}(Y)(P_{k}\otimes I_{H}))d^{-k} \\
=  \textrm{tr}_{{(\varphi}_{\infty}(M)\otimes I_{H})'}(\delta_{V}(Y)(\sum_{k=0}^{m}(P_{k}\otimes I_{H})d^{-k})).\label{eq:Popescu 1.4}\end{multline}
The passage from the first line to the second simply reflects the
properties of the trace and the fact that $\textrm{tr}_{{(\varphi}_{\infty}(M)\otimes I_{H})'}$
restricts to $\textrm{tr}_{{(\varphi}_{k}(M)\otimes I_{H})'}$ on
$P_{k}(\varphi_{\infty}(M)\otimes I_{H})'P_{k}$. Equation (\ref{eq:Popescu 1.4})
is an analogue of Popescu's equation (1.4) in \cite{gP01}. If $d\geq1$,
and if we replace $Y$ by $K(\eta)K(\eta)^{*}$ in equation (\ref{eq:Popescu 1.4}),
then we may take the limit as $m\to\infty$. The left hand side tends
to $\kappa(\eta)$ by Theorem \ref{CurvaturePoisson} (equation (\ref{eq:Curv2})),
while the right hand side tends to $\textrm{tr}_{(\varphi_{\infty}(M)\otimes I_{H})'}\{\delta_{V}[K(\eta)K(\eta)*]N\}$.

\end{proof}

\end{document}